\documentclass{siamltex}
\usepackage{amsmath,amssymb,amsfonts,latexsym,graphicx,color,multirow,footnote}

\newtheorem{algorithm}{Algorithm}

\title{A FEAST Algorithm with oblique projection for generalized eigenvalue problems }

\author{Guojian Yin\thanks{Department of Mathematics, The Chinese University of Hong Kong, Shatin, Hong  Kong ({\tt guojianyin@gmail.com}).} \and Raymond H. Chan\thanks{Department of Mathematics, The Chinese University of Hong Kong, Shatin, Hong Kong. Research is supported in part by HKRGC GRF Grant No. CUHK400412, HKRGC CRF Grant No. CUHK2/CRF/11G,
HKRGC AoE Grant AoE/M-05/12, CUHK DAG No. 4053007, and CUHK FIS Grant No. 1902036 ({\tt rchan@math.cuhk.edu.hk}).}  \and Man-Chung Yeung\thanks{Department of Mathematics, University of Wyoming,
Dept. 3036, 1000 East University Avenue, Laramie, WY 82071, USA ({\tt myeung@uwyo.edu}).}  }
\begin{document}
\maketitle

\begin{abstract}
The contour-integral based eigensolvers are the recent efforts for computing the eigenvalues inside a given region in the complex plane. The best-known members  are the Sakurai-Sugiura (SS) method, its stable version CIRR, and the FEAST algorithm. An attractive computational advantage of these methods is that they are easily parallelizable. The FEAST algorithm was developed for the generalized Hermitian eigenvalue problems. It is stable and accurate. However, it may fail when applied to non-Hermitian problems. In this paper, we extend the FEAST algorithm to non-Hermitian problems. The approach can be summarized as follows: (i) to construct a particular contour integral to form a subspace containing the desired eigenspace, and (ii) to use the oblique projection technique to extract desired eigenpairs with appropriately chosen test subspace. The related mathematical framework is established. We also address some implementation issues such as how to choose a suitable starting matrix and design good stopping criteria. Numerical experiments are
provided to illustrate that our method is stable and efficient.
\end{abstract}



\section{Introduction}
Consider the generalized eigenvalue problem
\begin{equation}\label{eq:1-1}
 A {\bf x} = \lambda  B {\bf x},
\end{equation}
where  $ A,  B \in {\mathbb C}^{n \times n}$. The scalars $\lambda \in \mathbb{C}$
and the associated vectors ${\bf x} \in {\mathbb C}^n, {\bf x}\neq 0$, are
called the eigenvalues and eigenvectors, respectively. In this paper,
we are concerned with computing the eigenvalues of (\ref{eq:1-1})
that are located inside a given region in the complex plane together with their eigenvectors.

Large-scale generalized eigenvalue problems arise in various areas of
science and engineering,
such as  dynamic analysis of structures \cite{GLS94}, determination of the
linearized stability of 3-D fluid flows {\cite{BS06}}, the electron energy
and position problems in quantum chemistry \cite{FH74}, the widely used
principal component analysis \cite{STL11}, and the linear discriminant
analysis in statistical data analysis \cite{CLM12}.
In some applications, it is not the whole spectrum but
rather a significant part of it is of interest to the users. For example, in the electronic structure calculations of materials \cite{SCS10},
it is required to compute the lowest portion of
the spectrum of (\ref{eq:1-1}); and
in the model reduction of a linear dynamical system, one only needs to know
the response over a range of frequencies, see \cite{BDDRV00, GGD96}.

Solving \eqref{eq:1-1} is a very challenging problem, even if
there are various practical methods and software available, see \cite{BDDRV00}.
When $A$ and $B$ have no special structures and the whole
spectrum is required, the QZ method \cite{MS73}
is the most widely used method. It uses a sequence of unitary equivalence
transformations to reduce the original pair $(A, B)$ to generalized Schur
form. The algorithm is numerically stable. However its computational cost is expensive, requiring about
$46n^3$ floating point operations \cite{gvl}.

There are several methods for computing only part of the spectrum of (\ref{eq:1-1}). The rational Krylov subspace method
approximates all eigenvalues in a union of regions around the chosen
shifts \cite{Ruhe98}. However, it needs locking, purging and implicit
restart techniques, which
create difficulties in practical implementation. The divide-and-conquer
approaches, which are based on the sign-function or inverse-free techniques, are
also popular choices \cite{BDG97}. However, these approaches always suffer
from slow convergence or poor stability \cite{NH13}.
When $A$ and $B$ are Hermitian, a shifted block
Lanczos algorithm was proposed  in \cite{GLS94} to compute the eigenvalues
contained in a given interval. The authors designed a shift strategy combining with the LDL decomposition  to
guarantee that all eigenvalues in the given interval can be found.
But  the shift strategy is complicated in practical implementation and its efficiency depends on the distribution of the spectrum.

The methods based on contour integral are recent efforts for solving partial spectrum of (\ref{eq:1-1}). When (\ref{eq:1-1}) is a diagonalizable and non-degenerate system,  a contour integral method, called the Sakurai-Sugiura (SS) method, was proposed in \cite{ss} for evaluating the eigenvalues inside a specific domain.  In this method, the original problem (\ref{eq:1-1}) is reduced to a small eigenproblem with Hankel matrices. However, since Hankel matrices are usually ill-conditioned \cite{BGL07}, the SS method always suffers from numerical instability \cite{AT14, ST07}. Later in \cite{ST07}, Sakurai {\it et al.} used the Rayleigh-Ritz procedure to replace the Hankel matrix approach to get a more stable algorithm called CIRR.
In \cite{ISN10} and \cite{IS10} the block versions of SS and CIRR were proposed respectively to make SS and CIRR also available for degenerate systems. It was shown that SS and CIRR, as well as their corresponding block variants, can be regarded as the Krylov subspace techniques \cite{IDS14, ISN10}.

Recently in \cite{polizzi}, Polizzi proposed another eigenproblem solver based on contour integral, called FEAST, to solve (\ref{eq:1-1})
under the assumptions that $A$ and $B$ are
Hermitian and $B$ is positive definite, i.e, (\ref{eq:1-1}) is a Hermitian system or $zB-A$ is a definite matrix pencil \cite{BDDRV00}. His algorithm
computes all eigenvalues inside a given interval, along
with their associated eigenvectors. The FEAST algorithm is accurate and reliable, see \cite{kramer} for more details. It was shown that the FEAST algorithm can be understood as a standard subspace iteration with the Rayleigh-Ritz procedure \cite{TP13}.

The FEAST algorithm originally was proposed for Hermitian problems. However, when it is applied to non-Hermitian problems, it may fail to find the desired eigenvalues. A simple example (Example 3.1) will be given later to illustrate this. Motivated by this fact, our goal is to generalize the FEAST algorithm to non-Hermitian problems and establish the related mathematical framework. The only requirement in our method is that the corresponding matrix pencil $zB-A$ is regular, i.e., ${\rm det}(zB-A)$ is not identically zero for all $z \in \mathbb{C}$. In other words, our new method can deal with the most common generalized eigenvalue problems \cite{BDDRV00}. Unlike FEAST which uses the Rayleigh-Ritz procedure to extract the desired eigenpairs, our generalized FEAST algorithm uses the oblique projection technique to find them.

One of the main drawbacks of the contour-integral based algorithms, including ours, is that the information about the number of desired eigenvalues has to be known a priori. This is because we need this information (i) to choose
an appropriate size for the starting matrix to start the methods, and (ii) to  determine whether all desired eigenvalues are captured when  the methods stop.
In this paper, we use a way similar to that proposed in \cite{SFT} to find an upper bound of the number of eigenvalues inside the target region. It will help to choose the starting matrix.  We also provide good stopping
criteria which not only can guarantee that all desired eigenpairs are captured,
but also can give the accuracy that our method can achieve. With these efforts, our method is applicable in practical implementation. Comparisons with {\sc Matlab}'s \texttt{eig} command and the block version of CIRR \cite{IS10, SFT} show that our method is an efficient and stable solver for large generalized eigenvalue problems.

The outline of the paper is as follows. In Section 2,
we briefly describe two typical contour-integral based eigensolvers: the CIRR method \cite{ST07} and the FEAST algorithm \cite{polizzi}. In Section 3, we extend the FEAST algorithm to non-Hermitian problems and establish the related mathematical framework.
In Section 4, we present a way to find a suitable upper bound of the number
of eigenvalues inside the prescribed region and give the stopping criteria. Then we present the complete algorithm of our method. In Section 5, numerical experiments are reported to illustrate the efficiency and stability of our method.

Throughout the paper, we use the following notation and terminology.
The subspace spanned by the columns of a matrix $ X$ is
denoted by ${\rm span}\{ X\}$. The rank and conjugate transpose of $X$ are 
denoted by $\rank(X)$ and $X^*$ respectively. We denote the submatrix
consisting of the first $i$ rows and the first $j$
columns of $X$ by $X_{(1:i,1:j)}$, the submatrixs consisting
of the first $j$ columns of $X$ and  the first $i$ rows of $X$ by $X_{(:,1:j)}$ and $X_{(1:i,:)}$ respectively. The algorithms are presented in \textsc{Matlab} style.

\section{Two typical contour-integral based eigensolvers}
In this section, we briefly introduce two typical contour-integral based eigensolvers: the CIRR method and the FEAST algorithm. Before starting our discussion, we present some facts about matrix pencil which will play important roles in the derivation of these methods.

Recall that the matrix pencil $(z B-A)$ is regular if ${\rm det}(zB-A)$ is not identically zero for all $z \in \mathbb{C}$. The Weierstrass canonical form for regular matrix pencils is given below.

\begin{theorem} [\cite{G59}] \label{thm2.1}
Let $zB-A$ be a regular matrix pencil of order $n$. Then there exist nonsingular matrices $S$ and $T \in \mathbb{C}^{n\times n}$ such that
\begin{equation}\label{eq:8-3-1}
TAS = \begin{bmatrix}
  J_d    & 0   \\
   0   & I_{n-d}
\end{bmatrix}  \quad {\rm and} \quad TBS= \begin{bmatrix}
   I_d   & 0   \\
  0    & N_{n-d}
\end{bmatrix},
\end{equation}
where $J_d$ is a $d\times d$ matrix in Jordan canonical form
with its diagonal entries corresponding to the eigenvalues of $zB-A$, $N_{n-d}$ is an $(n-d)\times (n-d)$ nilpotent matrix also in Jordan canonical form, and $I_d$ denotes the identity matrix of order $d$.
\end{theorem}

Let $J_d$ be of the form
\begin{equation}\label{equ:7-17-1}
J_d =   \left[ \begin{array}{cccc}
 J_{d_1}(\lambda_1) & 0 & \cdots & 0\\
0 &  J_{d_2}(\lambda_2) & \cdots & 0\\
\vdots & \vdots& \ddots &  \vdots\\
0 & 0 & \cdots &  J_{d_m}(\lambda_m)
\end{array}
\right]
\end{equation}
where $\sum_{i=1}^m d_i  =d$ and
$ J_{d_i}(\lambda_i)$ are $d_i \times d_i$ matrices of the form
$$\begin{array}{ccc}
 J_{d_i}(\lambda_i) = \left[ \begin{array}{ccccc}
\lambda_i & 1 & 0& \cdots  & 0\\
0 & \lambda_i & 1 & & \vdots\\
 & \ddots & \ddots & \ddots &0\\
\vdots & & \ddots &\ddots &  1\\
0 &  \cdots && 0&\lambda_i
\end{array}
\right], & & i = 1, 2, \ldots, m
\end{array}
$$
with $\lambda_i$ being the eigenvalues. Here the $\lambda_i$
are not necessarily distinct and can be repeated according to
their multiplicities.

Let $N_{(n-d)}$ be of the form
\begin{equation*}\label{equ:7-17-11}
N_{(n-d)} =   \left[ \begin{array}{cccc}
N_{d^{\prime}_1} & 0 & \cdots & 0\\
0 &  N_{d^{\prime}_2} & \cdots & 0\\
\vdots & \vdots& \ddots &  \vdots\\
0 & 0 & \cdots &  N_{d^{\prime}_{m^{\prime}}}
\end{array}
\right],
\end{equation*}
where $\sum_{i=1}^{m^{\prime}} d^{\prime}_i  =n-d$ and
$N_{d^{\prime}_i}$ are $d^{\prime}_i \times d^{\prime}_i$ matrices of the form
$$\begin{array}{ccc}
N_{d^{\prime}_i} = \left[ \begin{array}{ccccc}
0 & 1 & 0& \cdots  & 0\\
0 &0 & 1 & & \vdots\\
 & \ddots & \ddots & \ddots &0\\
\vdots & & \ddots &\ddots &  1\\
0 &  \cdots && 0&0
\end{array}
\right], & & i = 1, 2, \ldots, m^{\prime}.
\end{array}
$$

Let us partition $S$ into block form
\begin{equation}\label{eq:3-4-1}
S = [S_1, S_2, \ldots, S_ m, S_{m+1}],
\end{equation}
where $S_i\in \mathbb{C}^{n\times d_i}$, $1\le i\le  m$ and $S_{m+1}\in \mathbb{C}^{n\times (n-d)}$. One can easily verify by \eqref{eq:8-3-1}
(or see \eqref{Jordan})
that the first column in each $S_i$, $i = 1,\ldots, m$, is the eigenvector associated with the eigenvalue $\lambda_i$ of \eqref{eq:1-1}.

\subsection{The CIRR method}
In \cite{ss}, Sakurai {\it et al.} used a moment-based technique to formulate a contour-integral based method, which now is known as SS, for finding the eigenvalues of (\ref{eq:1-1}) inside a given region. Since the SS method suffers from numerical instability, a stable version, called the CIRR method, was later developed using the Rayleigh-Ritz procedure \cite{IS10, ST07}.

Below we show how to use the CIRR method to compute the eigenvalues inside $\Gamma$, which is a given positively
oriented simple closed curve in the complex plane.
Without loss of generality, let the set of eigenvalues of (\ref{eq:1-1}) enclosed by $\Gamma$ be $\{\lambda_1, \ldots, \lambda_l\}$, and $s: = d_1+d_2+\cdots + d_{l}$ be the number of eigenvalues inside $\Gamma$ with multiplicity taken into account.
Define the contour integrals
\begin{equation} \label{eq:2-2-21}
  F_i:= \dfrac{1}{2\pi \sqrt{-1}}\oint_{\Gamma}z^i (z B- A)^{-1} dz,\quad i=0,1, \ldots.
\end{equation}
It was shown in \cite{IS10} that
\begin{equation}\label{eq:2-19-1}
F_i = S_{(:,1:s)}(J_{(1:s, 1:s)})^iT_{(1:s,:)}
\end{equation}
where $S$ and $T$ are given by Theorem \ref{thm2.1}.
The CIRR method uses the Rayleigh-Ritz procedure to extract the eigenpairs inside $\Gamma$ \cite{ST07}. Originally, it was derived under the assumptions that (\ref{eq:1-1}) is a Hermitian system and the desired eigenvalues are distinct, i.e., they are non-degenerate \cite{ST07}. However, based on the following theorem,  the CIRR method was extended to non-Hermitian cases in  \cite{IS10}.

\vspace{2mm}
\begin{theorem}[\cite{IS10}]\label{Th:2-13-1}
Let $L$, $D \in \mathbb{C}^{n\times t}$, $t\ge s$, be arbitrary matrices, and $R = F_0D$, where $F_0$ is defined in (\ref{eq:2-2-21}). A projected matrix pencil $z\hat{B}-\hat{A}$ is defined by $\hat{B}=L^*BR$ and $\hat{A}=L^*AR$. If the ranks of both $L^*(T^{-1})_{(:,1:s)}$ and $T_{(1:s, :)}D$ are $s$, the non-singular part of the projected matrix pencil is equivalent to $zI_s-(J_{d})_{(1:s, 1:s)}$.
\end{theorem}

\vspace{2mm}
\noindent
Theorem \ref{Th:2-13-1} says that the desired eigenvalues $\{\lambda_i\}_{i=1}^s$ can be solved by computing the eigenvalues of the projected eigenproblem $z\hat{B}-\hat{A}$, if the ranks of both $L^*(T^{-1})_{(:,1:s)}$ and $T_{(1:s, :)}D$ are $s$.

Define the columns of $D$ to be
\begin{equation}\label{eq:5-30-1}
D_{(:, i)} = (T^{-1})_{(:, 1:s)}(J_{(1:s, 1:s)})^{i-1}T_{(1:s,:)}v,\quad i =1, 2, \ldots, t.
\end{equation}
It is shown in \cite{ISN10} that the rank of $T_{(1:s, :)}D$ is $s$, if all the elements of $T_{(1:s,:)}v$ are non-zero and there is no degeneracy in $J_{(1:s, 1:s)}$. By (\ref{eq:2-2-21}) and (\ref{eq:5-30-1}), we have
\begin{equation}\label{eq:3-10-1}
R_{(:, i)} = F_0D_{(:, i)} = F_iv, \quad i = 1,\ldots, t.
\end{equation}
Based on Theorem \ref{Th:2-13-1}, the right Ritz space is spanned by the vectors $\{F_iv\}_{i =0}^{t-1}$.  As for the left Ritz space, in the CIRR for non-Hermitian problems, it is chosen to be the same as the right one.  

In order to remove the restriction on the non-degeneracy in $J_{(1:s, 1:s)}$, a block CIRR method was also proposed in \cite{IS10}, where the random vector $v$ is replaced by a random matrix $Y\in \mathbb{C}^{n\times h}$ of appropriate dimension. The right and the left Ritz spaces are spanned by the vectors $\{F_iY\}_{i =0}^{g-1}$, where $g$ is a positive integer satisfying $hg \ge s$. Then all eigenvalues of $h^{\prime}$-order degeneracy, $h^{\prime}\le h$, can be found \cite{IS10, SFT}. Obviously, the main task of the block CIRR method is to evaluate $\{F_i Y\}_{i = 0}^{g-1}$. In practice, $F_i Y$ are computed approximately by a quadrature scheme according to (\ref{eq:2-2-21}).
Below is  the block CIRR algorithm for non-Hermitian problems.

\begin{algorithm}\label{alg:CIRR}
Input matrices $A$ and $B$,
a random matrix ${Y} \in \mathbb{C}^{n\times h}$, and a positive integer $g$ satisfying $t =hg \ge s$.
The function ``{\sc Block\_CIRR}" computes eigenpairs of
(\ref{eq:1-1}) that are located inside $\Gamma$,
and they are output in the vector $\Lambda_s$ and the matrix $X_s$.
\end{algorithm}
\vspace{.2cm}
\begin{tabbing}
x\=xxx\= xxx\=xxx\=xxx\=xxx\=xxx\kill
\> Function $[ \Lambda_s,  X_s]$ = {\sc Block\_CIRR}$( A,  B,   Y, g, \Gamma)$\\
\>1.\> Compute $R_i = F_{i} Y, i = 0, 1, \ldots, g-1, $ approximately by a quadrature scheme.\\
\>2.\> Compute the singular value decomposition: $[R_0,\ldots, R_{h-1}] = U\Sigma V^*$.\\
\>3.\> Set $\bar{A}=  U^* A  U$ and $\bar{B}=  U^* B  U$.\\
\>4.\> Solve the generalized eigenproblem of size $t$: $\bar{A} {{\bf y}}=\bar{\lambda}  \bar{B} {{\bf y}}$,
to obtain the \\
\>\>  eigenpairs $\{(\bar{\lambda}_i, {\bf y}_i)\}_{i = 1}^{t}$.\\
\>5.\> Compute $\bar{\bf {x}}_i=   U{{\bf y}}_i, i =1,2,\ldots t$, and select $s$ approximate eigenpairs.
\end{tabbing}
\vspace{.2cm}

We see that for Algorithm \ref{alg:CIRR}, in order to choose the parameters $h$ and $g$, one needs to know  $s$ and the degrees of degeneracy of the desired eigenvalues. In the recent article \cite{SFT}, Sakurai {\it et al.} gave a method to choose a suitable  $h$ for fixed $g$. We will describe the method in Section \ref{sec:number}. Also in \cite{SFT}, the authors suggested to perform iterative refinement in case the eigenpairs computed by Algorithm \ref{alg:CIRR}
cannot attain the prescribed accuracy.

\subsection{The FEAST algorithm}

In this section, we give a brief introduction to the FEAST algorithm
{\cite{polizzi}} due to Polizzi. The FEAST algorithm was formulated under the assumptions $A$ and $B$ are Hermitian and $B$ is positive definite, in which case the eigenvalues of (\ref{eq:1-1}) are real-valued \cite{BDDRV00}.
It is used to find all eigenvalues of
(\ref{eq:1-1}) within a specified interval, say $[\sigma_1,\sigma_2]$,
and their associated eigenvectors.
Here we also assume that the desired eigenvalues are $\lambda_1\le \lambda_2\le\ldots \le\lambda_s$. Let $\Gamma$ be any contour that contains $\{\lambda_i\}_{i=1}^s$ inside. For example, $\Gamma$ can be the circle with center at $c = (\sigma_1+\sigma_2)/2$ and radius $r = (\sigma_2-\sigma_1)/2$.

When $zB-A$ is a definite matrix pencil, the Weierstrass canonical form  (\ref{eq:8-3-1})  is reduced to
\begin{equation}\label{eq:2-27-1}
T A S= \Lambda ={\rm diag}([\lambda_1,\lambda_2
\cdots, \lambda_n])\quad {\rm and}  \quad  T B S=  I_n.
\end{equation}
It is easy to see that $T = S^*$, and the columns $S_{(:, i)}$ are the eigenvectors corresponding to $\lambda_i$, $i = 1,\ldots, n$. Therefore
according to (\ref{eq:2-19-1}), we have
\begin{equation}\label{eq:2-27-2}
F_0 = S_{(:,1:s)}T_{(1:s,:)}= S_{(:,1:s)}( S_{(:,1:s)})^{*}.
\end{equation}

Let $ W: =  F_0 Y = S_{(:,1:s)}( S_{(:,1:s)})^{*}Y$, where $Y$ is an $n\times t$ matrix with $t\ge s$. Then $ W$ forms a basis for the desired eigenspace ${\rm span}\{S_{(:,1:s)}\}$ if $( S_{(:,1:s)})^{*}Y$ is  full-rank. In the FEAST algorithm, the elements of $Y$ are chosen to be random numbers to increase
the chance that $ W$ may form a basis for ${\rm span}\{  S_{(:,1:s)}\}$.
According to the Rayleigh-Ritz procedure \cite{polizzi, stewart},
the problem (\ref{eq:1-1}) is transformed to the problem of computing the eigenpairs of the smaller generalized eigenvalue problem
$$\hat{A} {\bf y} = \lambda \hat{B} {\bf y}$$
of size $t\times t$, where $\hat{A}= W^* A  W$, $\hat{B}=  W^* B  W$, and $t< n$.

Now the task is to get the basis $W= F_0Y$. Since $S_{(:,1:s)}$ is unknown,
one cannot use (\ref{eq:2-27-2}).  Instead $W$ is computed by (\ref{eq:2-2-21})  numerically using a quadrature scheme such as the Gauss-Legendre quadrature rule \cite{DR84}.
The complete FEAST algorithm is given as follows.

\begin{algorithm}\label{alg:8}
Input Hermitian matrices $A$ and $B$ with $B$ being positive definite,
a random matrix $Y \in \mathbb{R}^{n\times t}$,
where $t\geq s$,
the circle $\Gamma$ enclosing the interval $[\sigma_1,\sigma_2]$, and
a convergence  tolerance $\epsilon$.
The function ``{\sc Feast}" computes eigenpairs $(\hat{\lambda}_i, \hat{\bf x}_i)$ of
(\ref{eq:1-1}) that satisfy
\begin{equation}\label{con-cre}
\hat{\lambda}_i \in [\sigma_1,\sigma_2]
\quad {\rm and} \quad \sum_{i = 1}^s \hat{\lambda}_i< \epsilon,
\end{equation}
and they are output in the vector $\Lambda_s$ and the matrix $X_s$.
\end{algorithm}
\vspace{.2cm}
\begin{tabbing}
x\=xxx\= xxx\=xxx\=xxx\=xxx\=xxx\kill
\> Function $[ \Lambda_s,  X_s]$ = {\sc Feast}$( A,  B,  Y, \Gamma, \epsilon)$\\
\>1.\> Compute $ {W}=F_0Y $ approximately by the Gauss-Legendre quadrature rule.\\
\>2.\> Set $\hat{A}=  W^* A  W$ and $\hat{B}=  W^* B  W$.\\
\>3.\> Solve the generalized eigenproblem of size $t$: $\hat{A} {{\bf y}}=\hat{\lambda}  \hat{B} {{\bf y}}$,
to obtain the \\
\>\>  eigenpairs $\{(\hat{\lambda}_i, {\bf y}_i)\}_{i = 1}^{t}$.\\
\>4.\> Compute $\hat{\bf {x}}_i=   W{{\bf y}}_i, i =1,2,\ldots t$.\\
\>5.\> Check if $\{(\hat{\lambda}_i, \hat{\bf x}_i)\}_{i=1}^t$ satisfy the convergence criteria (\ref{con-cre}). If $s$ eigenpairs \\
\>\> satisfy (\ref{con-cre}), stop. Otherwise, set $X_t=  [\hat{\bf x}_1,\hat{\bf x}_2,\ldots,\hat{\bf x}_{t}]$ and $Y=BX_t$,  then\\
\>\> go back to Step 1.
\end{tabbing}
\vspace{.2cm}

The FEAST algorithm
can be understood as a standard subspace iteration combining with the Rayleigh-Ritz procedure \cite{TP13}. It is an accurate and reliable technique {\cite{kramer}}.
However, like CIRR, in practice we have to know $s$ in advance in order to choose $t$ for the starting matrix $Y$ and to determine whether all desired eigenpairs are found. In \cite{TP13}, a technique was presented to find an estimation of $s$ under the conditions that $A$ and $B$ are Hermitian and $B$ is positive definite.
The estimation can help to select an appropriate starting matrix $Y$.

\section{A FEAST method with oblique projection}\label{sec:feast} The FEAST method was originally developed for generalized Hermitian eigenvalue problems. We note that it may fail when applied to non-Hermitian problems, as is shown by the following example.

\vspace{2mm}\noindent
\textsc{Example 3.1}: Let $A$ and $B$ be defined as follows:
$$
 A=
\begin{pmatrix}
 0 &  0  & 0  & 5   \\
   0 & 0  & 2  & 0  \\
   0 & 0.5  & 0  & 0  \\
   0.2 & 0  & 0  & 0
\end{pmatrix}
\quad \mbox{and}\quad  B=\begin{pmatrix}
    0  & 0  & 0  & 1  \\
    0  & 0  & 1  & 0  \\
    0  & 1  & 0  & 0  \\
    1  & 0  & 0  & 0
\end{pmatrix}.
$$
Then (\ref{eq:8-3-1}) holds with
$$
T=\begin{pmatrix}
    1  & 0  & 0  & 0  \\
    0  & 1  & 0  & 0  \\
    0  & 0  & 1  & 0  \\
    0  & 0  & 0  & 1
\end{pmatrix},\quad \mbox{and}\quad  S=\begin{pmatrix}
    0  & 0  & 0  & 1  \\
    0  & 0  & 1  & 0  \\
    0  & 1  & 0  & 0  \\
    1  & 0  & 0  & 0
\end{pmatrix}.
$$
In fact, $TAS = {\rm diag}([5, 2, 0.5, 0.2])$ and $TBS = I_4$.
Suppose we want to find the eigenvalues of $A{\bf x}=\lambda B{\bf x}$ lying inside the unit circle. Obviously, the eigenvalues of interest are $0.2$ and $0.5$, and the corresponding
eigenvectors are $S_{(:, 4)}=[1, 0, 0, 0]^*$ and $S_{(:, 3)}=[0, 1, 0, 0]^*$. By (\ref{eq:2-19-1}) and (\ref{eq:2-27-2}), it is easy to check that all projected matrices $\hat{A}$ and $\hat{B}$ in the FEAST algorithm are zero for any given random matrix $Y$. Hence any complex number will be an eigenvalue of the corresponding projected eigenproblem $\hat{A}{\bf y} = \lambda \hat{B}{\bf y}$!

In view of the above example, we now extend the FEAST algorithm to non-Hermitian problems. The only requirement of our method is that the corresponding matrix pencil $ z B-A$ is regular. Hence our method can deal with the most common generalized eigenvalue problems \cite{BDDRV00}.  Below we establish the related mathematical framework.

Again without loss of generality, we let the eigenvalues of (\ref{eq:1-1}) enclosed by $\Gamma$ be $\{\lambda_1, \ldots, \lambda_l\}$, and $s: = d_1+d_2+\cdots + d_{l}$ be the number of eigenvalues inside $\Gamma$ with multiplicity taken into account. Define the contour integral
\begin{equation}\label{equ:7-26-8}
 Q: = F_0 = \frac{1}{2\pi \sqrt{-1}} \oint_\Gamma  (zB-A)^{-1}B  dz.
 \end{equation}
For $z \ne \lambda_i$, according to (\ref{equ:7-17-1})
the matrix $zI_d-J_d$ is invertible. Hence
by (\ref{eq:8-3-1}), the resolvent operator is given by
\begin{eqnarray}\label{eq:8-3-4}
(zB-A)^{-1}B & = & S\begin{bmatrix}
  (zI_d-J_d)^{-1}    &   0 \\
    0  &  (zN_{n-d}-I_{n-d})^{-1}
\end{bmatrix} TB \nonumber \\
 & = & S\begin{bmatrix}
  (zI_d-J_d)^{-1}    &   0 \\
    0  &  (zN_{n-d}-I_{n-d})^{-1}
\end{bmatrix} \begin{bmatrix}
  I_d    &   0 \\
    0  &  N_{n-d}
\end{bmatrix}S^{-1}\nonumber \\
 & = & S\begin{bmatrix}
  (zI_d-J_d)^{-1}    &   0 \\
    0  &  (zN_{n-d}-I_{n-d})^{-1}N_{n-d}
\end{bmatrix}S^{-1}.
\end{eqnarray}
Notice that the first diagonal block $(zI_d-J_d)^{-1}$
in (\ref{eq:8-3-4}) is of the block diagonal form
where each diagonal sub-block is of the form:
\begin{equation}\label{equ:7-26-3}
(z  I_{d_i} -  J_{d_i}(\lambda_i))^{-1} = \left[ \begin{array}{cccc}
 \displaystyle{\frac{1}{z-\lambda_i}}& \displaystyle{\frac{1}{(z-\lambda_i)^2}}& \cdots &
 \displaystyle{\frac{1}{(z-\lambda_i)^{d_i}}}\\
0 & \displaystyle{\frac{1}{z-\lambda_i}} & \cdots & 
\displaystyle{\frac{1}{(z-\lambda_i)^{d_i-1}}}\\
\vdots& \vdots& \ddots & 
\vdots\\
0 & 0 & \cdots 
&\displaystyle{\frac{1}{z-\lambda_i}}\\
\end{array}
\right].
\end{equation}
Similarly, the second diagonal block in (\ref{eq:8-3-4}) is also of the block diagonal form where each diagonal sub-block is of the form:
\begin{equation}\label{eq:7-3-7}
(z N_{d^{\prime}_i}-I_{d^{\prime}_i})^{-1}N_{d^{\prime}_i}=\begin{bmatrix}
   0   &  -1 & -z & \cdots & -z^{d^{\prime}_i-2}  \\
    0  & 0 & -1 &\cdots & -z^{d^{\prime}_i-3}\\
    0  & 0 & 0 &\cdots &-z^{d^{\prime}_i-4}\\
    \vdots  & \vdots & \vdots &\ddots & \vdots\\
    0 & 0 & 0 &\cdots &0
\end{bmatrix}.
\end{equation}

Then, according to the residue theorem in complex analysis \cite{Rudin}, it follows from (\ref{eq:8-3-4})--(\ref{eq:7-3-7}) that
\begin{equation}\label{equ:7-26-4}
 Q = \frac{1}{2\pi \sqrt{-1}} \oint_\Gamma  (zB-A)^{-1}B  dz =  S\left[ \begin{array}{cc}
  I_s & 0\\
0 & 0
 \end{array} \right]  S^{-1} =  S_{(:,1:s)} ( S^{-1})_{(1:s,:)}.
\end{equation}
Using the remark following (\ref{eq:3-4-1}), we know that ${\cal K}:={\rm span}\{S_{(:,1:s)}\}$ contains
the eigenspace corresponding to the eigenvalues $\{\lambda_1, \ldots, \lambda_l\}$.
Since $ Q^2 =  Q$,  $Q$ is a spectral projector onto ${\cal K}$.
Define $U:=QY$, where $Y$ is an
appropriately chosen matrix so that $U$ forms a basis for ${\cal K}$.
As in FEAST, we choose $Y$ randomly, and we
show below that the resulting $U$ does form a basis for ${\cal K}$.

\vspace{2mm}
\begin{lemma}\label{lem:8-8-1}
Let $Y\in\mathbb{R}^{n\times s}$. If the entries of
$Y$ are random numbers from a continuous distribution and that they are independent and identically distributed (i.i.d.), then with probability 1, the matrix
$( S^{-1})_{(1:s,:)}Y$ is nonsingular.
\end{lemma}

\vspace{2mm}
\begin{proof}
Let $Z= ( S^{-1})_{(1:s,:)} Y$. Consider $|\det(Z)|^2$, the square of the absolute value of $\det(Z)$, as a real coefficient polynomial in the elements of $Y$. We now show that the polynomial is
non-zero, i.e. $|\det(Z)|^2  \not\equiv 0$.

Since the rank of $( S^{-1})_{(1:s,:)}$ is $s$, it has an $s\times s$ nonsingular submatrix. Without loss of generality, let the left $s \times s$
submatrix of $( S^{-1})_{(1:s,:)}$ be nonsingular. Then set $Y=[ I_s,  0]^T$, we have
$\det(Z) = \det(( S^{-1})_{(1:s,:)}Y) \ne 0$. Therefore the polynomial $|\det(Z)|^2$ is not identically zero.
Hence the set ${\cal Z}$ of zeros of the polynomial $|\det(Z)|^2 $ is
of measure zero in ${\mathbb R}^{ns}$ according to \cite[Prop.~4]{joubert2}. When $Y$ is randomly picked, it is with probability 1 that $Y \not\in {\cal Z} $, or equivalently, $|\det(( S^{-1})_{(1:s,:)} Y)|^2  \ne 0$. \qquad
\end{proof}

Since
\begin{equation}\label{eq:6-12}
U =QY=  S_{(:,1:s)} ( S^{-1})_{(1:s,:)} Y,
\end{equation}
and $( S^{-1})_{(1:s,:)} Y$ is nonsingular by {Lemma} \ref{lem:8-8-1},
the columns of $U$ form a basis for $\mathcal{K}$. Our next step is to  project the original problem (\ref{eq:1-1}) onto a small subspace where we can
extract the required eigenpairs.
Unlike CIRR and FEAST which use the Rayleigh-Ritz procedure to extract the desired eigenpairs, here we resort to the oblique projection method, namely, the Petrov-Galerkin condition \cite{BDDRV00, saad}. Since $\mathcal{K}$ contains the eigenspace corresponding to the eigenvalues inside $\Gamma$, it is natural to choose $\mathcal{K}$ as the search subspace. Our next task is to seek an appropriate test subspace.

Let us further partition each $S_i$ in (\ref{eq:3-4-1}) into $S_i=[{\bf s}_1^i, {\bf s}_2^i, \ldots ,{\bf s}_{d_i}^i]$ with
${\bf s}_j^i \in {\mathbb C}^{n}$, $1\le i \le m, 1\le j \le d_i$.
Notice that by (\ref{eq:8-3-1}), we have for any eigenvalue $\lambda_i$,
$1\le i \le m$,
\begin{eqnarray}
(\lambda_i B-A)S \begin{bmatrix}
  I_d &   0 \\
    0  &  N_{n-d}
\end{bmatrix}
& = &T^{-1}\begin{bmatrix}
  \lambda_i I_d-J_d   &   0 \\
    0  &  (\lambda_iN_{n-d}-I_{n-d})N_{n-d}
\end{bmatrix} \nonumber \\
&=&
T^{-1}\begin{bmatrix}
  I_d  &   0 \\
    0  &  N_{n-d}
\end{bmatrix}\begin{bmatrix}
  \lambda_i I_d-J_d   &   0 \\
    0  &  \lambda_iN_{n-d}-I_{n-d}
\end{bmatrix} \nonumber \\
&=&
BS\begin{bmatrix}
  \lambda_i I_d-J_d   &   0 \\
    0  &  \lambda_iN_{n-d}-I_{n-d}
\end{bmatrix}. \label{eq:6-1-1}
\end{eqnarray}
By comparing the first $d$ columns on both sides above, we get
\begin{equation}\label{Jordan}
(\lambda_i B-A){\bf s}_j^{i} = B{\bf s}_{j-1}^{i}, \quad 1\le j \le d_i,
\quad 1 \le i \le m,
\end{equation}
with ${\bf s}_0^{i} \equiv {\bf 0}$.
In particular, ${\bf s}_1^{i}$ is the eigenvector corresponding
to the eigenvalue $\lambda_i$ for all $1 \le i \le m$.
From (\ref{Jordan}), we see that
$A {\cal K} \subseteq B {\cal K}$. Therefore we choose $B {\cal K}$ as the test subspace. The Petrov-Galerkin condition then becomes:
\begin{equation}\label{eq:6-3-11}
(A{\bf x}_i-{\lambda}_iB{\bf x}_i)\  \bot \ B\mathcal{K}, \quad
1 \le i \le l,
\end{equation}
with ${\lambda}_i\in \mathbb{C}$ and ${\bf x}_i\in \mathcal{K}$.

Now we are in the position to find a basis for $B\mathcal{K}$. From (\ref{eq:8-3-1}), we know that  the rank of $BS_{(:,1:s)}$ is $s$.  Hence by {Lemma} \ref{lem:8-8-1}, $BU=BS_{(:,1:s)} ( S^{-1})_{(1:s,:)} Y$ is full-rank, which implies that $BU$ forms a basis for $B\mathcal{K}$.
Recall that $U$ forms a basis for $\mathcal{K}$, and we seek an
${\bf x}_i\in\mathcal{K}$. Therefore (\ref{eq:6-3-11}) can be written in matrix form
\begin{equation}\label{eq:6-3-12}
(BU)^*(AU{\bf y}_i-{\lambda}_iBU{\bf y}_i)=0,
\end{equation}
where ${\bf y}_i \in \mathbb{C}^s$ satisfying ${\bf x}_i=U{\bf y}_i$. 
Accordingly, we get the projected eigenproblem
\begin{equation}\label{eq:6-3-13}
\tilde{ A}{\bf y}=\tilde{\lambda}\tilde{ B}{\bf y},
\end{equation}
with
\begin{equation}\label{equ:7-10}
\tilde{ A}=  (BU)^* A U\quad \mbox{and}\quad \tilde{ B}=  (BU)^* B U.
\end{equation}
Our method is to compute the desired eigenpairs of \eqref{eq:1-1} by solving the projected eigenproblem (\ref{eq:6-3-13}). The theory behind our method is
given in the next theorem.

\vspace{2mm}
\begin{theorem}\label{Th:8-3-3}
\begin{enumerate}
\item[(a)] Let $\{(\tilde{\lambda}_i, {\bf y}_i)\}_{i=1}^s$ be $s$ eigenpairs of
the projected eigenproblem (\ref{eq:6-3-13}). Then
$\{(\tilde{\lambda}_i, U{\bf y}_i)\}_{i=1}^s$ are the eigenpairs of (\ref{eq:1-1})  located inside ${\Gamma}$.
\item[(b)] If $\mathcal{Y}_{\tilde{\lambda}_i}$ is the eigenspace of (\ref{eq:6-3-13}) corresponding to the eigenvalue $\tilde{\lambda}_i$, then $U \mathcal{Y}_{\tilde{\lambda}_i}$ is the eigenspace of (\ref{eq:1-1}) corresponding to the eigenvalue $\tilde{\lambda}_i$.
\end{enumerate}
\end{theorem}

\vspace{2mm}
\begin{proof}
(a): 
First, since $A \mathcal{K} \subseteq B \mathcal{K}$ and $BU$ forms a basis for $B \mathcal{K}$, there exist vectors ${\bf q}_i \in \mathbb{C}^{s}, 1\le i \le s,$ such that
\begin{equation*}\label{eq:6-3-14}
AU{\bf y}_i-\tilde{\lambda}_i BU{\bf y}_i= BU {\bf q}_i,
\quad 1 \le i \le s.
\end{equation*}
By (\ref{eq:6-3-12}), we have $(BU)^*(BU){\bf q}_i=0$.
Since $BU$ is full-rank, ${\bf q}_i=0$. Consequently,
$AU{\bf y}_i=\tilde{\lambda}_i BU{\bf y}_i$.
Thus $\{(\tilde{\lambda}_i, U{\bf y}_i)\}_{i=1}^s$  are the eigenpairs of (\ref{eq:1-1}).

Next we want to show that
$\{\tilde{\lambda}_i\}_{i=1}^s$ are exactly the $s$
eigenvalues of (\ref{eq:1-1}) inside ${\Gamma}$.  By (\ref{eq:6-1-1}), we can easily verify
that
\begin{equation*}\label{eq:6-3-16}
(z B- A)S_{(:,1:s)}=BS_{(:,1:s)}(zI_s-(J_d)_{(1:s,1:s)}), \quad \forall z \in {\mathbb C}.
\end{equation*}
Hence by the definitions in (\ref{equ:7-10}) and (\ref{eq:6-12}), we have
\begin{eqnarray}\label{eq:4-2-18}
    z\tilde{B}-\tilde{A} &= & (BU)^*(zB-A)  S_{(:,1:s)} ( S^{-1})_{(1:s,:)}  Y \nonumber \\
      & = & (BU)^*BS_{(:,1:s)}(zI_s-(J_d)_{(1:s,1:s)})( S^{-1})_{(1:s,:)}  Y.
\end{eqnarray}
Therefore,
\begin{equation*}\label{eq:6-3-17}
 {\rm det}(z\tilde{B}-\tilde{A}) ={\rm det}((BU)^*BS_{(:,1:s)}) {\rm det}(zI_s-(J_d)_{(1:s,1:s)}){\rm det}(( S^{-1})_{(1:s,:)}  Y).
\end{equation*}
Since $BS_{(:,1:s)}$ is of full rank $s$ (see (\ref{eq:8-3-1})), and $( S^{-1})_{(1:s,:)}  Y$ is nonsingular by {Lemma} \ref{lem:8-8-1},
the matrix
$(BU)^*BS_{(:,1:s)}  = ((S^{-1})_{(1:s,:)}Y)^* (BS_{(:,1:s)})^* BS_{(:,1:s)}$ is nonsingular. Hence
${\rm det}(z\tilde{B}-\tilde{A})=0$ if and only if
${\rm det}(zI_s-(J_d)_{(1:s,1:s)})=0$. By the special structure of
$ (J_d)_{(1:s,1:s)}$ (see (\ref{equ:7-17-1})), the zeros of the determinant ${\rm det}(z\tilde{B}-\tilde{A})$ are precisely
$\{\lambda_i\}_{i=1}^l$ with multiplicities $\{d_i\}_{i=1}^l$ respectively. Therefore, $\{\tilde{\lambda}_i \}_{i=1}^s$ are precisely all the eigenvalues of (\ref{eq:1-1}) inside ${\Gamma}$.

(b):
Let $\mathcal{X}_{\tilde{\lambda}_i}$ be the eigenspace of (1.1) corresponding to the eigenvalue $\tilde{\lambda}_i$. Then $U \mathcal{Y}_{\tilde{\lambda}_i} \subseteq \mathcal{X}_{\tilde{\lambda}_i}$ by Part(a).
From (\ref{equ:7-17-1}) and (\ref{eq:4-2-18}),
it can be seen that $\dim(\mathcal{Y}_{\tilde{\lambda}_i})$ is equal to the number of Jordan blocks in $(J_d)_{(1:s,1:s)}$ corresponding to the eigenvalue $\tilde{\lambda}_i $. On the other hand, the later coincides with the number of Jordan blocks in $J_d$ corresponding to $\tilde{\lambda}_i$, which is equal to $\dim(\mathcal{X}_{\tilde{\lambda}_i})$. Therefore,  we have
$\dim(\mathcal{Y}_{\tilde{\lambda}_i}) = \dim(\mathcal{X}_{\tilde{\lambda}_i})$.
Since $U$ has full column rank, it follows that
$$\dim(U \mathcal{Y}_{\tilde{\lambda}_i}) = \dim(\mathcal{Y}_{\tilde{\lambda}_i}) = \dim(\mathcal{X}_{\tilde{\lambda}_i}).$$
Therefore, we have $U \mathcal{Y}_{\tilde{\lambda}_i} = \mathcal{X}_{\tilde{\lambda}_i}$.
\end{proof}

Thus computing the eigenpairs of (\ref{eq:1-1}) inside ${\Gamma}$ is transformed into
computing the eigenpairs of the small $s \times s$ projected problem (\ref{eq:6-3-13}), which can be solved by standard solvers in LAPACK \cite{BDDRV00, demmel}, such as \texttt{xGGES} and \texttt{xGGEV} \cite{Lapack}.

In order to construct the projected eigenproblem (\ref{eq:6-3-13}), the most important task is to compute the matrix $U$  in (\ref{eq:6-12}). In practice, we have to compute $U$ by using the contour integral in (\ref{equ:7-26-4}), i.e.
\begin{equation}\label{eq:2-13}
 U =QY=\frac{1}{2\pi \sqrt{-1}} \oint_\Gamma (zB-A)^{-1}BdzY,
\end{equation}
which can be approximated by using for example the Gauss-Legendre quadrature rule.

We summarize our above derivation into the following algorithm.\\

\begin{algorithm}\label{alg:5}
Input $ A,  B \in {\mathbb C}^{n \times n}$, an i.i.d. random
matrix $Y \in {\mathbb R}^{n \times  t}$
where $ t\geq s$, a closed curve $\Gamma$, a convergence tolerance
$\epsilon$, and ``${\texttt{{max\_iter}}}$" to control the
maximum number of iterations. The function ``{\sc Eigenpairs}"
computes eigenpairs $(\tilde{\lambda}_i, \tilde{{\bf x}}_i)$ of (\ref{eq:1-1}) that satisfies
\begin{equation}\label{con-cre_2}
\tilde{\lambda}_i \ {\rm inside} \ \Gamma
\quad {\rm and} \quad \frac{\| A\tilde{{\bf x}}_i - \tilde{\lambda}_i  B\tilde{{\bf x}}_i\|_2}{\| A \tilde{{\bf x}}_i\|_2+
\| B \tilde{{\bf x}}_i\|_2} < \epsilon.
\end{equation}
The results are stored in the vector $\Lambda$ and the matrix $X$. \end{algorithm}
\vspace{.2cm}
\begin{tabbing}
x\=xxx\= xxx\=xxx\=xxx\=xxx\=xxx\kill
\> Function $[ \Lambda,  X] = \textsc{Eigenpairs}( A,  B,  Y, \Gamma, \epsilon, \texttt{max\_iter})$\\
\>1.\>For $k = 1,\cdots, \texttt{max\_iter}$\\
\>2.\>\> Compute $U$ in (\ref{eq:2-13}) approximately by the Gauss-Legendre quadrature rule. \\
\>3.\>\> Compute QR decompositions: $U = U_1R_1$ and $BU = U_2R_2.$\\
\>4.\>\>Form $ \tilde{A} =  U_2^*  A  U_1$ and $\tilde{B} =  U_2^*  B  U_1$.\\
\>5.\>\> Solve the projected eigenproblem $\tilde{A} {\bf y} = \tilde{\lambda} \tilde{B} {\bf y}$ of size $ t$ to obtain eigenpairs\\
\>\>\>$\{(\tilde{\lambda}_i, {\bf y}_i)\}_{i=1}^{t}$. Set  $ \tilde{{\bf x}}_i =  U_1{\bf y}_i, i = 1, 2, \ldots,  t$. \\
\>6.\>\> Set $ \Lambda=\left[ \  \right]$ and $X=\left[ \  \right]$.\\
\>7.\>\> For $i = 1: t$\\
\>8.\>\>\>If $(\tilde{\lambda}_i,\tilde{ {\bf x}}_i)$ satisfies (\ref{con-cre_2}), then $ \Lambda = [\Lambda, \tilde{\lambda}_i]$ and $ X =[X, \tilde{{\bf x}}_i]$.\\
\>9.\>\>End\\
\>10.\>\>If there are $s$ eigenpairs satisfy (\ref{con-cre_2}), stop. Otherwise, set $Y = U_1$.\\
\>11.\>End. \\
\end{tabbing}

{Algorithm} \ref{alg:5} faces the same issue occurred in CIRR and FEAST, that is we have to know $s$ in advance in order to choose $t$ for the starting matrix $Y$ and to determine whether all desired eigenpairs are found. More precisely, the number of columns $t$ of  $Y$ should satisfy $t\geq s$. This is because if $t < s$, then $\rank( U)  <s$. Consequently, the columns of $ U$ cannot form a basis for  $\mathcal{K}$. Moreover, since $s$ is unknown a priori, it is also hard to decide whether all desired eigenvalues are found and therefore it is hard to decide when to stop the algorithm. In the next section, we present strategies to address these two problems, which will make the resulting algorithm applicable to practical implementation.

\section{Our Algorithm} In this section, we first introduce a
method to find an upper bound for $s$. Our method is similar to a technique proposed in \cite{SFT}. Next we design stopping criteria to guarantee all eigenvalues are captured. After that, we present the complete algorithm.

\subsection{Finding an Upper Bound for the Number of Eigenvalues inside $\Gamma$}\label{sec:number}
In the following, by ``$M\sim {\sf N}_{p\times q}(0,1)$", we mean
$M$ is a $p\times q$ matrix with
i.i.d. entries drawn from the standard normal distribution ${\sf N}(0,1)$.

In \cite{futa}, an approach was proposed for finding an estimation of $s$. Here we derive a similar method. Let  $Y \sim {\sf N}_{n \times p}(0,1)$. One can easily verify that the mean
${\mathbb E} [{\rm trace}(Y^* M Y)] = p \cdot {\rm trace}(M)$ for any $n\times n$ matrix $M$.
In particular, by (\ref{equ:7-26-4}) and (\ref{eq:6-12}),
\begin{eqnarray}
\frac1p{\mathbb E}[{\rm trace} (Y^*U) ]&=&\frac1p{\mathbb E}[{\rm trace} (Y^*  Q Y)] \nonumber \\
& = &{\rm trace}( Q) =
{\rm trace}( S_{(:,1:s)} ( S^{-1})_{(1:s,:)})\nonumber \\
 &=&{\rm trace}(( S^{-1})_{(1:s,:)}  S_{(:,1:s)}) = {\rm trace}( I_s)=s . \label{trace}
\end{eqnarray}
So $s_0:=\frac1p{\rm trace} (Y^*U)$ is a good initial estimation  of $s$. In \cite{futa}, the entries of $Y$ are taken to be $1$ or $-1$ with equal probability.

In Algorithm \ref{alg:5}, we need to choose an upper bound $t$ of $s$
for the starting matrix $Y \in {\mathbb R}^{n \times t}$. However, in practice, $s_0$ may be less than $s$, and we may not know this fact as we do not know  $s$. Below we present a way to find a better estimate $t$ based on $s_0$.

Recall that in CIRR, for a given integer $g$, we need to select an $h$ such that $t = hg \ge s$.  In \cite{SFT}, a method was proposed for selecting a suitable $h$. The method works as follows. Assume an estimation $s_0$ of $s$ is available. Let $Y\in \mathbb{C}^{n\times h}$, where $h = \lceil \frac{\kappa s_0}{g}\rceil$ and $\kappa >1$. Compute $P_h=[F_0BY, F_1BY,\ldots, F_{g-1}BY]$ (see (\ref{eq:2-2-21})) and the minimum singular value $\sigma_{\min}$ of $P_h$. If  $\sigma_{\min}$  is not small, then $h$ is increased until $\sigma_{\min}$ of the updated $P_h$ is small enough. Similar to this idea, we determine an upper bound
$t$ for our Algorithm \ref{alg:5} by using the numerical rank of $U = QY$,  where $Y \sim {\sf N}_{n \times t}(0,1)$. The rationale behind our method is as follows.

Let $s^{\dagger}$ be a positive integer and $ Y_{s^{\dagger}} \sim {\sf N}_{n\times s^{\dagger}}(0,1)$. Consider
$$
U_{s^{\dagger}}=QY_{s^{\dagger}} = S_{(:,1:s)} (S^{-1})_{(1:s,:)} Y_{s^*}= \frac{1}{2\pi \sqrt{-1}} \oint_\Gamma  (zB-A)^{-1}B  Y_{s^{\dagger}} dz.
$$
Then $ U_{s^{\dagger}} \in {\mathbb C}^{n \times s^\star}$ is the projection of $ Y_{s^{\dagger}}$ onto $\mathcal{K}$, and consequently, $\rank( U_{s^{\dagger}} ) \leq s$. With this in mind, if $\rank( U_{s^{\dagger}} ) = s^{\dagger}$, it obviously means that $s^{\dagger} \leq s$; and
we will increase $s^\star$ and repeat the process. Otherwise, if
$\rank( U_{s^{\dagger}} ) < s^{\dagger}$, we can conclude that $s= \rank( U_{s^{\dagger}} )$ with the help of {Lemma} \ref{lem:8-8-1}; and thereby $s < s^{\dagger}$.
Below we give the algorithm for finding $t$.\\

\begin{algorithm}\label{alg:2}
Input  an increasing factor $\alpha > 1$ and the size $p$ of sample vectors. The function ``\textsc{Search}" outputs $ t$ (an upper bound of $s$) and the projection matrix $ U_1  \in {\mathbb C}^{n \times  t}$
onto $\mathcal{K}$.
\end{algorithm}
\vspace{.2cm}
\begin{tabbing}
x\=xxx\= xxx\=xxx\=xxx\=xxx\=xxx\kill
\>Function $[ U_1, t] = \textsc{Search}( A,  B, \Gamma,  \alpha, p)$ \\
\>1. \> Pick $ Y_0 \sim {\sf N}_{n \times p}(0,1)$
 and compute $\displaystyle{ U = \frac{1}{2\pi \sqrt{-1}} \oint_\Gamma  (zB-A)^{-1}B  Y_0 dz}$\\
 \>\> approximately by the Gauss-Legendre quadrature rule .\\
\>2.\>  Set $s_0 = \lceil \frac{1}{p}{\rm trace}( Y_0^*  U)\rceil$ and $s^{\dagger} = \min(\max(p, s_0), n)$.\\
\>3.\> If $s^{\dagger} > p$\\
\>4.\>\> Pick $\tilde{Y} \sim {\sf N}_{n \times (s^{\dagger} -p)}(0,1)$ and compute
$\displaystyle{\tilde{U} = \frac{1}{2\pi \sqrt{-1}} \oint_\Gamma (zB-A)^{-1}B \hat{Y} dz}$ \\
\>\>\> approximately the Gauss-Legendre quadrature rule.\\
\>5.\>\> Augment $V$ to $U$ to form $ U = [ U, \tilde{U} ] \in {\mathbb C}^{n \times s^{\dagger}}$.\\
\>6.\>Else\\
\>7.\>\>Set $s^{\dagger} = p$. \\
\>8.\> End\\
\>9.\> Compute $ U=U_1R\Pi$: the rank-revealing QR decomposition \cite{tony} of $ U$. \\
\>10.\> Set  $ t= {\rm rank}( R)$. If $ t<s^{\dagger}$,  stop.  Otherwise, set $p = t,  U = U_1$ and  \\
\>\>$s^{\dagger} = \lceil \alpha t\rceil$. Then go to Step 3.
\end{tabbing}
\vspace{.2cm}

In Algorithm
\ref{alg:2}, we use the rank-revealing QR decomposition method \cite{tony, gvl} to detect whether $U$ is numerically rank deficient (line 9). If this case occurs, it means that the subspace spanned by $U_1$  already contains $\mathcal{K}$, so we stop the procedure and get a good upper bound $t$. We remark that {Algorithm} \ref{alg:2} can be treated as the first iteration of {Algorithm} \ref{alg:5} because $U_1$ from {Algorithm} \ref{alg:2}
is a projection onto $\mathcal{K}$, therefore it can be taken as the $U$ in line 2 of {Algorithm} \ref{alg:5}.

\subsection{The Stopping Criteria} \label{sect4.2}
Although by {Algorithm} \ref{alg:2} we can obtain an upper bound $t$ for $s$, the actual value of $s$ is still unknown. Thus {Algorithm} \ref{alg:5} is
still impractical because it is hard to determine whether all desired eigenpairs are captured. Below we present simple but efficient stopping criteria to guarantee this. It
can also give the accuracy that the algorithm have achieved.

In {Algorithm} \ref{alg:5}, there are $t$ eigenvalues being solved in each iteration, and hence $s$ of them are the eigenvalues we sought
and $(t-s)$ of them are spurious eigenvalues. Those spurious eigenvalues outside $\Gamma$ can easily be detected by checking the values of
their coordinates. It is only the spurious eigenvalues that are inside $\Gamma$
that we need special attention. As the iteration
progresses, the accuracy of the $s$ desired eigenvalues will improve steadily while the
accuracy of the spurious eigenvalues will not. Therefore after some iterations, there will be a gap
in the accuracy between the desired eigenvalues and the spurious eigenvalues. Based on this observation, we choose a tolerance $\eta$ for detecting the number of desired eigenvalues inside $\Gamma$. If after some iterations, there are $s^{\prime}$ eigenvalues inside $\Gamma$ whose accuracy are smaller than $\eta$, and the value of $s^{\prime}$ is unchanged in
two consecutive iterations, then we set $s= s^{\prime}$. In our experiments, we set $\eta =10^{-3}$ and we can determine $s$ in one or two iterations.

After determining $s$, we will continue with the iterations
so as to improve the accuracy
of the $s$ desired eigenvalues. The algorithm will stop when all
$s$ eigenvalues meet the user-prescribed tolerance $\epsilon$, see (\ref{con-cre_2}).
To avoid the situation where $\epsilon$ is set too small for the given
problem, we also stop the iterations when (i) the iteration
number reaches a prescribed maximum $\mathrm{\texttt{max\_iter}}$, or (ii) the overall accuracy of the $s$ eigenvalues is not improved
from one iteration to the next. Since we are monitoring the
accuracy of the $s$ desired eigenvalues in each iteration, we
have an estimate of how accurate the $s$ eigenvalues are
when the algorithm stops.

More detailed demonstration about the idea behind the stopping criteria will be given in Section 5.

\subsection{The Complete Algorithm}
In this section, we  present the complete algorithm which is based on {Algorithm} \ref{alg:5}, Algorithm \ref{alg:2}, and the above stopping criteria.
Then we discuss some implementation issues pertaining to the algorithm.

\begin{algorithm}\label{alg:6}
Input tolerance $\eta$ for detecting the spurious eigenvalues, and the
 tolerance $\epsilon$ for the accuracy of the eigenpairs. The function ``\textsc{GFEAST}" computes all the eigenvalues $\lambda$
of (\ref{eq:1-1}) inside $\Gamma$ and their associated eigenvectors ${\bf x}$ .
The computed $\lambda$ and ${\bf x}$  are stored in vector $ \Lambda$ and matrix $ X$ respectively. The $\mathrm{\texttt{flag}}$ is set to $1$ if there are $s$ eigenpairs $(\lambda, {\bf x})$ satisfying (\ref{con-cre_2}),  $0$ if the overall accuracy is not improved from
 one iteration to the next; $-1$ if the maximum number of iterations $\mathrm{\texttt{max\_iter}}$ is reached.
\end{algorithm}
\vspace{.2cm}
\begin{tabbing}
x\=xxx\= xxx\=xxx\=xxx\=xxx\=xxx\kill
\>Function $[ \Lambda,  X, \texttt{Err}, \texttt{flag}] = \textsc{GFEAST}( A,  B, \Gamma,  \alpha, p, \epsilon, \eta, \mathrm{\texttt{max\_iter}})$ \\
\>1.\> Call $[ U_1, t] = \textsc{Search}( A,  B, \Gamma,  \alpha, p)$ to obtain an upper bound $ t$ of the exact\\
\>\>  number $s$ of the eigenvalues inside $\Gamma$, and a projection $U_1$ onto $\mathcal{K}$.\\
\>2.\>Compute the QR decomposition: $BU_1=U_2 R_2$. \\
\>3.\> Set $e(0) = 0 $ and $c(0)=0$. \\
\>4.\> For $k = 1, 2, \cdots, \texttt{max\_iter}$ \\
\>5.\>\> Form $\tilde{A}=U^*_2AU_1$ and $\tilde{B}=U^*_2BU_1$.\\
\>6.\>\> Solve the projected eigenproblem $\tilde{A} {\bf y} = \lambda \tilde{B} {\bf y}$ of size $ t$ to obtain eigenpairs\\
\>\>\>$\{(\lambda_i, {\bf y}_i)\}_{i=1}^{t}$. Set  ${\bf x}_i =  U_1{\bf y}_i, i = 1, 2, \ldots,  t$. \\
\>7.\>\> Set $r =\left[ \  \right] , \Lambda^{(k)}=\left[ \  \right]$, $X^{(k)}=\left[ \  \right]$ and $c(k) = 0$.\\
\>8.\>\> For $i = 1: t$\\
\>9.\>\>\> Compute $r _i= \| A{\bf x}_i - \lambda_i  B{\bf x}_i\|_2/(\| A {\bf x}_i\|_2+\| B {\bf x}_i\|_2).$\\
\>10.\>\>\> If $\lambda_i$ inside $\Gamma$ and $ r_i < \eta$, then $c(k) = c(k)+1$, $r = [r, r_i]$,  \\
\>\>\>\>$X^{(k)} = [X^{(k)}, {\bf x}_i]$  and $\Lambda^{(k)} =[\Lambda^{(k)}, \lambda_i]$.\\
\>11.\>\>End\\
\>12.\>\> Set $e(k) = \max(r)$.\\
\>13.\>\>If $c(k)=c(k-1)$ and $e(k) < \epsilon$, output $ \Lambda= \Lambda^{(k)}$ and $ X= X^{(k)}$, \\
\>\>\>$\texttt{Err} = e(k), \texttt{flag} =1$. Stop.\\
\>14.\>\>If $c(k)=c(k-1)$ and $e(k) > e(k-1)$, output $ \Lambda= \Lambda^{(k-1)}$  \\
\>\>\>and $ X= X^{(k-1)}$, $\texttt{Err} = e(k-1), \texttt{flag} =0$. Stop.\\
 \>15.\>\> If $k = \texttt{max\_iter}$, output $ \Lambda= \Lambda^{(k)}$ and $ X= X^{(k)}$,  \\
\>\>\>$\texttt{Err} = e(k), \texttt{flag} =-1$. Stop.\\
\>16.\> \>Compute  $U = \frac{1}{2\pi \sqrt{-1}} \oint_\Gamma  (zB-A)^{-1}B dzU_1$ approximately.\\
\>17.\>\> Compute QR decompositions: $U = U_1R_1$ and $BU = U_2R_2.$\\
\>18.\>End
\end{tabbing}
\vspace{.2cm}

Below we give some remarks on Algorithm \ref{alg:6}.
\begin{enumerate}
\item Since the columns of the matrix $U_1$ obtained from the function \textsc{Search}  are orthonormalized, we only need to compute the QR decomposition for $BU_1$ in line 2.
\item In line 10, we keep only those eigenpairs inside $\Gamma$ whose accuracy are less than $\eta$,
  and we consider them to be the $s$ desired eigenpairs  of (\ref{eq:1-1}) inside $\Gamma$.
\item
Lines 13 to 15 are the three stopping criteria.  In line 13,
we stop when $e(k)$, the overall accuracy of all $s$ desired eigenvalues,  is less than
$\epsilon$. In line 14, we stop when $e(k)$ is not
improved from the $(k-1)$th iteration to the $k$th iteration. In line 15, we stop when the maximum number of iterations is reached.
 \end{enumerate}

In each iteration, the dominant work is to compute the projection $U$. In the case when $\Gamma$ is a circle with center $\gamma$ and radius $\rho$, we can
compute the contour integral in line 16 by using the $q$-point Gauss-Legendre quadrature on $\Gamma$. More precisely
\begin{equation}\label{eq:8-4}
  U = \dfrac{1}{2\pi \sqrt{-1}}\oint_{\Gamma}  (z B- A)^{-1}Bdz U_1
 \approx \frac{1}{2} \sum^{q}_{j=1}\omega_j(z_j-\gamma)(z_j  B- A)^{-1}B  U_1,
\end{equation}
where $z_j=\gamma+\rho e^{\sqrt{-1}\theta_j}$, $\theta_j=(1+ t_j)\pi$, and
$t_j$ is the $j$th Gaussian node with associated weight $\omega_j$.
 Accordingly, it requires us to solve $q$ generalized shifted linear systems of the form
\begin{equation}\label{eq:4-4}
(z_j  B- A)\widehat{ U}_j =  BU_1\quad z_j\in \mathbb{C}, \quad 1 \le j \le q.
\end{equation}
When $\Gamma$ is an irregular closed curve, we can choose a circle $\widetilde{\Gamma}$ such that $\widetilde{\Gamma}$ encloses $\Gamma$. We compute all eigenpairs inside $\widetilde{\Gamma}$ and
then determine the eigenpairs that are indeed inside $\Gamma$.

Like other contour-integral based methods, our algorithm replaces the difficulty of solving the eigenvalue problem (\ref{eq:1-1}) by the difficulty of solving the linear systems (\ref{eq:4-4}). One has considerable freedom to choose different approaches to solve (\ref{eq:4-4}) based on the properties of the matrices in (\ref{eq:4-4}), such as the Krylov subspace based methods \cite{MB10, SS07, SHZF12}. Since (\ref{eq:4-4}) are generalized shifted systems with multiple right-hand sides,  the direct  methods, such as the sparse Gaussian LU factorization, are also highly recommended. Notice that once we obtain the LU factors, they can be reused when we solve  (\ref{eq:4-4}) in the subsequent iterations. Moreover, since the quadrature nodes $z_j, j =1,\ldots , q,$ are independent, and
 the columns of the right-hand sides are also independent, our algorithm has a good potential to be parallelized.

\section{Numerical Experiments}\label{sec:experiments}
In this section, we give some numerical experiments to illustrate the efficiency
of our  method for computing the eigenpairs of (\ref{eq:1-1}) inside
a given contour $\Gamma$.  All computations are carried out in \textsc{Matlab}
version R2012b on a MacBook with an Intel Core i5 2.5 GHz processor and 8 GB RAM.
The test matrices are from the Matrix Market collection \cite{DGL89}. They are
the real-world problems from scientific and
engineering applications. In fact, Problem 2 serves as a classic testbed for generalized non-Hermitian eigenproblem \cite{BDDRV00}. The eigenvalues of Problem 3 occur in pairs, and hence the problem is well-known as a difficult eigenproblem \cite{DGL89}.
The matrix $B$ of Problem 4 is singular and Problems 6 and 7 are ill-conditioned.
The descriptions of these matrices are presented in
\textsc{Table} \ref{Tab:5-1}, where nnz denotes the number of non-zero entries
and \texttt{cond} denotes their condition numbers which are computed by \textsc{Matlab} function \texttt{condest}.

\begin{table}
\caption{Test problems from Matrix Market that are used in our experiments.}
\footnotesize{
\noindent
\begin{tabular}{c|llllc}
No.&Matrix & Size&nnz & Property & \texttt{condest}
\\ \hline
1& $ A$: BFW398A & $398$& $3678$&  unsymmetric &$7.58\times 10^{3}$\\
& $ B$: BFW398B  & $398$& $2910$& symmetric indefinite &$3.64\times 10^{1}$\\ \hline
2& $ A$: BFW782A& $782$ & $7514$& unsymmetric &$4.63\times 10^{3}$\\
& $ B$: BFW782B & $782$& $5982$&  symmetric indefinite &$3.05\times 10^{1}$\\ \hline
3&$ A$: PLAT1919 & $1919$& 17159&  symmetric indefinite & $1.40\times 10^{16}$\\
&$ B$: PLSK1919 & $1919$& 4831& skew symmetric & $1.07\times 10^{18}$\\ \hline
4&$ A$: BCSSTK13& $2003$& 42943& symmetric positive definite &$4.57\times 10^{10}$\\
&$ B$: BCSSTM13 & $2003$& 11973& symmetric positive semi-definite &Inf\\ \hline
5&$ A$: BCSSTK27& $1224$& 28675& symmetric positive definite & $7.71\times 10^{4}$\\
&$ B$: BCSSTM27 & $1224$& 28675& symmetric indefinite & $1.14\times 10^{10}$\\ \hline
6 & $ A$: MHD3200A & $3200$ & 68026 &unsymmetric & $2.02\times 10^{44}$\\
 & $ B$: MHD3200B & $3200$ & 18316& symmetric indefinite & $2.02\times 10^{13}$\\ \hline
7&$ A$: MHD4800A & $4800$& 102252& unsymmetric & $2.54\times 10^{57}$\\
& $ B$: MHD4800B & $4800$&  27520& symmetric indefinite &$1.03\times 10^{14}$
\end{tabular}}
\label{Tab:5-1}
\end{table}

In the numerical comparisons, we assume that the eigenvalues and eigenvectors
computed by the \textsc{Matlab} function \texttt{eig} in dense format are the accurate ones. We use Gauss-Legendre quadrature rule \cite{DR84} with $q=16$ quadrature points on $\Gamma$ to compute the contour integrals \eqref{eq:8-4}.
As for solving the generalized shifted linear systems of the form (\ref{eq:4-4}),
we first use the \textsc{Matlab} function \texttt{lu} to compute the LU
decomposition of $ A-z_j B, j =1 ,2,\ldots, q$, and then perform the triangular
substitutions to get the corresponding solutions.

\textsf{Experiment\ 5.1 (Finding an upper bound for $s$):}
As stressed in the introduction and in Section 3, the information about the number of eigenvalues $s$
inside $\Gamma$ is crucial to the success of contour-integral based methods---they
all need an upper bound of $s$ to start the program with. Our {Algorithm} \ref{alg:2} is devoted to finding an upper bound $t$.
In this experiment, we test how accurate the computed $t$ are.
In the algorithm, the size $p$ of sample vectors is set as $50$, and the increasing factor $\alpha$
is chosen to be $1.5$. \textsc {Table} \ref{Tab:5-2} presents the results
for the test problems in \textsc {Table} \ref{Tab:5-1}.
In \textsc {Table} \ref{Tab:5-2}, the parameters $\gamma$ and $\rho$ denote the
center and the radius of the circle $\Gamma$ of each test problem respectively,
 $s$ denotes the true number of eigenvalues inside $\Gamma$ as obtained from \textsc{Matlab}
(by computing all eigenvalues and selecting those inside $\Gamma$),
 $s_0$ is the initial estimate obtained by the trace formula
(see (\ref{trace}) or line 2  in  {Algorithm} \ref{alg:2}), and $t$ is the upper bound that our {Algorithm} \ref{alg:2} gives.

From \textsc {Table} \ref{Tab:5-2}, we see that for the first five well-conditioned
problems, the estimates $s_0$ are good approximations of $s$ though it can underestimate
$s$ as in Problem 1. However for the ill-conditioned Problems 6 and 7,  $s_0$ are not good---it underestimates and overestimates $s$ by large margins.
However our {Algorithm}
\ref{alg:2} gives quite reasonable upper bounds $t$ in all seven problems.

\begin{table}
\centering
\caption{Results of {Algorithm} \ref{alg:2}: $s$ is the exact number of eigenvalues inside $\Gamma$
as computed by {\sc Matlab}, $s_0$ is the estimate of $s$ by using the trace formula (\ref{trace}) and $t$ is the upper bound computed by our {Algorithm} \ref{alg:2}.}
\footnotesize{
\noindent
\begin{tabular}{c|cc|ccc}
No.&$\gamma$&$\rho$&$s$& $s_0$&$t$  \\ \hline
1& $-5.0\times 10^5$&$2.0\times 10^5$  &123 &122& 137   \\
2&  $-6.0\times 10^{5}$&$3.0\times 10^{5}$  & 230&231& 262 \\
3&$0$ & $1.0\times 10^{-3}$&270 &277 &328 \\
4& $0$& $6.0\times 10^{5}$& 172 & 173 & 183  \\
5& $5.0\times 10^3$& $2.0\times 10^3$ & 107& 107& 118 \\
6 & $-4.0\times 10^1$ & $3.0\times 10^1$ & 162 &118& 178  \\
7& $-6.0$& $3.0$ & 169&3667& 186
\end{tabular}}
\label{Tab:5-2}
\end{table}

$\textsf{Experiment\ 5.2 (Separation of true and spurious eigenvalues):}$
Once we get an upper bound $t$ of the number of eigenvalues
inside $\Gamma$, we can separate the $(t-s)$ spurious eigenvalues
from the $s$ desired ones by using a threshold $\eta$ as 
explained in Section \ref{sect4.2}. In fact, the accuracy of the true
eigenvalues will continue to improve as the iteration progresses 
while those of the spurious eigenvalues will not. Therefore
the number of eigenvalues with accuracy better than $\eta$ will improve monotonically with each iteration, but once the number becomes constant 
in two consecutive iterations, we treat that number as our computed $s$. The accuracy of each eigenpair $(\lambda_i, {\bf x}_i)$ is measured by
\begin{equation}\label{eq:5-1}
r_i = \frac{\| A{\bf x}_i - \lambda_i  B{\bf x}_i\|_2}{\| A {\bf x}_i\|_2+\| B {\bf x}_i\|_2},
\quad 1 \le i \le t.
\end{equation}

\textsc {Table} \ref{Tab:5-3} gives the number of iterations required to
get our computed $s$ with the given $\eta$'s. We remark that our computed $s$ are
exactly the same as the true $s$ computed by \textsc{Matlab} so we do not write
them out again in \textsc {Table} \ref{Tab:5-3}. We see that if $\eta$ is set too small (e.g. $10^{-9}$), not all desired eigenvalues can attain such accuracy and we may not be able to get the true $s$. However, for $\eta=10^{-3}$, we need only one or two iterations to get the true $s$. Hence in the following experiments, we set $\eta=10^{-3}$.

\begin{table}
\centering
\caption{Number of iterations required to get the correct $s$ for different $\eta$'s.}
\footnotesize{
\noindent
\begin{tabular}{c|ccccccc}
$\eta$ & No.1  & No.2  & No.3 & No.4 & No.5 & No.6 & No.7 \\
\hline
$10^{-3}$ & $1$& $2$&$1$&$2$ & $2$ & $2$& $2$  \\
$10^{-5}$ & $2$ & $2$&$2$&$2$ &$3$& $3$& $3$  \\
$10^{-7}$ & $2$& $3$&$3$&$3$ &$3$& $4$ & $3$ \\
$10^{-9}$ & $3$& $3$&$-$&$-$ &$4$& $-$ & $-$
\end{tabular}}
\label{Tab:5-3}
\end{table}

\begin{figure}
\begin{center}
\includegraphics[width=13cm]{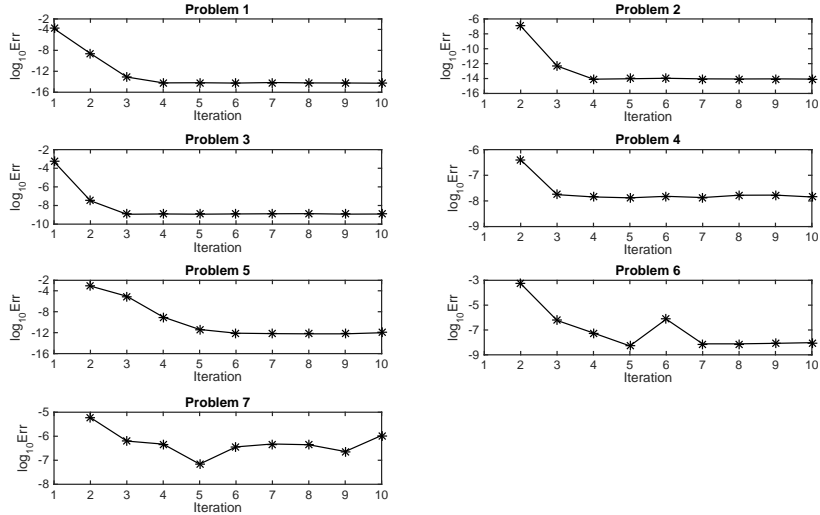}
\caption{The maxmum relative residual norms in different iterations.}
\label{Fig:5-3}
\end{center}
\end{figure}

$\textsf{Experiment\ 5.3 (Stopping criteria):}$
Here we illustrate the convergence behavior of our algorithm
and explain the stopping criteria we used. In \textsc{Fig.}~\ref{Fig:5-3}, 
we plot the maximum error of the $s$ desired eigenvalues 
$\texttt{Err}:=\max_{1\le i\le s} r_i$ for all seven problems,
starting from the iterations where $s$ is first determined until the
10th iteration (the iteration numbers where $s$ is first determined are
given in the first row of \textsc{Table} \ref{Tab:5-3}). We see from \textsc{Fig.}~\ref{Fig:5-3}
that $\texttt{Err}$ decreases monotonically and dramatically in the first few iterations
for all test problems.
Then it maintains at almost the same level for the first five problems
while it rebounds for the ill-conditioned problems Nos. 6 and 7.

Thus the stopping criteria in our  algorithm are
\begin{enumerate}
\item[(i)] when $\texttt{Err}$ in the current iteration is less than
a given tolerance $\epsilon$ (line 13 in {Algorithm}  \ref{alg:6}),
\item[(ii)] when $\texttt{Err}$ starts to increase again from one iteration to
the next (line 14 in {Algorithm}  \ref{alg:6}),
or
\item[(iii)] when the maximum number of iterations (line 15 in {Algorithm}  \ref{alg:6}) is reached.
\end{enumerate}
For example, if we set $\epsilon=10^{-10}$, the algorithm will stop in the sixth iteration
for both Problems 6 and 7, and return the eigenpairs it found in the 5th iteration.

\textsf{Experiment\ 5.4 (Comparisons with other methods):}
Here we compare our  method with two other methods both in terms of accuracy and
timing. 

We first compare our method  with  \textsc{Matlab}  function \texttt{eig}. We set $\epsilon=10^{-16}$ and  $\texttt{max\_iter}  = 10$ for our method. The goal is to examine the accuracy that our method can achieve. For all test problems, our 
method stops in Line 14 before reaching $\texttt{max\_iter}$.
The comparison of these two solvers are listed in \textsc {Table} \ref{Tab:5-5}. It is clear that our  algorithm can achieve higher accuracy when compared
with \texttt{eig} and can find all eigenvalues inside the target regions.
In terms of CPU time, except for Problem 1 where the size is small, our  algorithm runs significantly faster---though we should stress that \texttt{eig}
has to compute all eigenvalues while our method computes only those inside $\Gamma$.
\begin{table}
\centering
\caption{Comparison of   \texttt{eig} and our method.}
\footnotesize{
\noindent
\begin{tabular}{c|c|ccc|ccc}
\multirow{2}{*}{No.}  & \multirow{2}{*}{Size} & \multicolumn{3}{c|}{\texttt{eig}} &  \multicolumn{3}{c}{our method} \\ 
 & & $\texttt{Err}$ & \#eig &  Time (sec.) & $\texttt{Err}$ & \#eig &  Time (sec.) \\
\hline
1 & 398&$6.75\times 10^{-14}$&$123$&$ 0.67 $ & $6.02\times10^{-15}$ & 123 & 1.21 \\
2 & 782&$3.83\times 10^{-14}$& $230$ &$10.22$ & $8.70\times10^{-15}$ & $230$ & $6.16$ \\
3 &1919&$4.77\times 10^{-9}$& $270$ &$116.12$ & $1.07\times10^{-9}$ & $270$ & $31.24$ \\
4 &2003&$2.52\times 10^{-6}$& $172$ &$269.56$ & $1.21\times10^{-8}$ & $172$ & $29.10$ \\
5 &1224&$1.25\times 10^{-11}$& $107$ &$48.83$ & $4.86\times10^{-13}$ & $107$ & $8.39$ \\
6 &3200&$5.76\times 10^{-7}$& $162$ &$561.29$ & $5.22\times10^{-9}$ & $162$ & $24.90$ \\
7 &4800&$1.44\times 10^{-7}$& $169$ &$2725.43$ & $5.10\times10^{-8}$ & $169$ & $101.91$
\end{tabular}}
\label{Tab:5-5}
\end{table}

Finally we compare our method with the block version of the CIRR method ({\sc Block\_CIRR}), i.e., Algorithm \ref{alg:CIRR}. For the sake of fairness, in the test we apply the recently developed iterative refinement approach \cite{SFT} to {\sc Block\_CIRR}.
We set $\epsilon= 10^{-8}$ and $\texttt{max\_iter}=10$ for both our method and {\sc Block\_CIRR}.  The numerical results are reported in \textsc {Table} \ref{Tab:5-6}. We see that the {\sc Block\_CIRR} algorithm fails
for Problems 4, 6 and 7 where we recall that the matrix $B$ of Problem 4 is singular, and Problems 6 and 7 are ill-conditioned. Therefore,
we see that our algorithm is more accurate and stable when compared to the {\sc Block\_CIRR}.

However, {\sc Block\_CIRR} outperforms our method in terms of timing. The dominant computational cost in each iteration of both methods are the solution of $q=16$ linear systems of the form in (\ref{eq:4-4}).
But the number of right-hand sides in the {\sc Block\_CIRR} is always set to be a small number ($16$ in the tests here), while in our method, it is
$t$ $(\ge s)$. Consequently, the block-CIRR method always requires
less CPU time than our method. It is our future project to extend our
method to block form so as to minimize the number of right hand sides.

\begin{table}
\centering
\caption{Comparison of  {\sc Block\_CIRR} and our method.}
\footnotesize{
\noindent
\begin{tabular}{c|c|ccc|ccc}
\multirow{2}{*}{No.}&\multirow{2}{*}{$s$}    & \multicolumn{3}{c|}{{\sc Block\_CIRR}} &  \multicolumn{3}{c}{our method} \\
 &  & $\texttt{Err}$ & \#eig &  Time (sec.) & $\texttt{Err}$ & \#eig &  Time (sec.)\\
\hline
1 & $123$& $1.41\times10^{-9}$ & 123 & 0.68&$9.16\times 10^{-14}$&$123$&$ 0.72 $  \\
2 & $230$ & $9.47\times10^{-9}$ & $230$ & $5.14$&$1.70\times 10^{-11}$& $230$ &$3.23$ \\
3 &$270$ & $2.14\times10^{-9}$ & $270$ & $3.66$&$6.45\times 10^{-9}$& $270$ &$17.58$ \\
4 &$172$& $2.22\times10^{-2}$ & $172$ & $21.13$& $1.21\times 10^{-8}$& $172$ &$29.10$  \\
5&$107$ & $2.32\times10^{-12}$ & $107$ & $1.23$ &$6.03\times 10^{-13}$& $107$ &$4.21$\\
6 &$162$ & $7.23\times10^{-1}$ & $121$ & $4.84$ &$7.74\times 10^{-9}$& $162$ &$21.10$\\
7 &$169$ & $1.99\times10^{-1}$ & $120$ & $98.26$&$5.10\times 10^{-8}$& $169$ &$101.91$
\end{tabular}}
\label{Tab:5-6}
\end{table}

\section{Conclusions}  We develop a contour-integral based method which extends the FEAST algorithm to non-Hermitian problems. It can compute the eigenvalues lying inside a given region in the complex plane and their associated eigenvectors. To extract the desired eigenpairs, we use the oblique projection technique with appropriately chosen test subspace rather than the Rayleigh-Ritz procedure. The numerical experiments illustrate that our algorithm is fast and can achieve high accuracy. We also provide a way to find an upper bound of the number of eigenvalues
inside the contour, and give stopping criteria to guarantee that all eigenvalues are captured when the method stops. Our  algorithm is easily parallelizable.
How to further improve its numerical performance, and to extend it to nonlinear eigenproblems will be our future work.

\section{Acknowledgment}
We would like to thank Dr. Peter P. T. Tang who introduced the FEAST algorithm to us. We would also like to thank Professor Tetsuya Sakurai for many fruitful discussions and providing us the codes of the {\sc Block\_CIRR} method.

\end{document}